\documentclass{elsart}

\usepackage{amsmath,amssymb,amscd,graphicx,subfigure}

\newcommand{\bi}{\begin{itemize}}
\newcommand{\ei}{\end{itemize}}
\newcommand{\be}{\begin{enumerate}}
\newcommand{\ee}{\end{enumerate}}
\newcommand{\bc}{\begin{center}}
\newcommand{\ec}{\end{center}}
\newcommand{\bt}{\begin{tabular}}
\newcommand{\et}{\end{tabular}}
\newcommand{\ba}{\begin{array}}
\newcommand{\ea}{\end{array}}

\newcommand{\cs}{context-sensitive}
\newcommand{\csl}{context-sensitive language}
\newcommand{\cf}{context-free}
\newcommand{\lan}{language}
\newcommand{\lba}{linear-bounded automaton}
\newcommand{\ra}{\rightarrow}

\newcommand{\N}{\mathbb N}
\newcommand{\perm}{permutation}

\newcommand{\ALR}{ALR}
\newcommand{\Bona}{MR2078910}
\newcommand{\HU}{MR645539}
\newcommand{\GrigIndex}{MR1678812}
\newcommand{\ChomS}{MR0152391}
\newcommand{\regev}{regev}
\newcommand{\flajolet}{flajolet}
\newcommand{\Gessel}{MR1041448}
\newcommand{\Imm}{MR961049}
\newcommand{\Schut}{MR0135680}
\newcommand{\Stanley}{MR1676282}

\begin{document}

\begin{frontmatter}

\title{Pattern
avoiding permutations\\ are context-sensitive}

\author {Murray Elder\thanksref{grant}} 
\address{School~of
Mathematics and Statistics,\\ University of St~Andrews, St~Andrews,  Scotland } 

\thanks[grant]{Supported by EPSRC grant GR/S53503/01 }

\begin{abstract}
We prove that a variant of the insertion encoding of Albert, Linton and Ru\v{s}kuc 
for any class of pattern avoiding permutations is context-senstive.
It follows that every finitely based class of permutations bijects to a \csl.
\end{abstract}


\end{frontmatter}


\maketitle

\section{Introduction}
In this article we consider a connection between the fields of formal language theory from
computer science, and restricted permutations from combinatorics.           
In particular, we set out to describe the set of all permuations that 
avoid some pattern $q$ in terms of formal languages, in an attempt to
enumerate them. We prove that there is a bijection between the set of 
permutations which avoid $q$ and a context-sensitive language. This language 
is a variation of the ``insertion encoding'' of Albert, Linton and Ru\v{s}kuc \cite{\ALR}.

In formal language theory there is a heirachy of languages of increasing
complexity, starting at regular languages, then context-free, stack, indexed,
context-sensitive, decidable and recursively enumerable languages, each 
one strictly contained in the previous one.

In his paper on P-recursiveness and D-finiteness from 91 \cite{\Gessel} Gessel states
that ``another possible candidate for P-recursiveness is the problem 
of counting permutations with forbidden subsequences'', which B{\'o}na strengthens 
to a conjecture in \cite{\Bona}.
Sch{\"u}tzenberger proved that regular languages have rational
generating functions \cite{\Schut}, 
and  Chomsky and Sch{\"u}tzenberger proved that 
unambiguous context-free languages have algebraic generating functions \cite{\ChomS}. 
Beyond this, little is known 
about generating functions for the higher languages. A natural question is to   
ask whether there is a formal language correpsponding to D-finite generating
functions. We give examples of  context-sensitive languages       
whose generating functions are not D-finite.  

In spite of this, the connections between languages and pattern avoidance   
seem to be worth pursuing. An interesting general question is to find the 
language class of lowest complexity (in the heirachy or otherwise) such that 
the set of permatations avoiding any pattern can be described
by languages in this class. Since the generating function for 
$1234$-avoiding permutations 
is known to be non-algebraic (but D-finite)  (see \cite{\Bona} p. 210 and below), 
we need to look above the 
class of unambiguous context-free for an answer.

The article is organised as follows. In Section \ref{sec:start} we
 define the notion of pattern avoidance and basis for permutations, and 
then define \cs\ grammars and their machine counterparts. 
In  Section \ref{sec:encoding} we
 define the encoding, and in Section \ref{sec:main} we prove the main theorem, 
that the encoding is a
 \csl.

\section{Preliminaries}\label{sec:start}

A {\em permutation} is a finite ordered list of distinct real numbers. 
Two permutations $p_1p_2\ldots p_n$ and $q_1q_2\ldots q_n$ of the same 
length are said to be {\em order isomorphic} if $p_i<p_j$ if an only if
 $q_i<q_j$ for all $1\leq i<j\leq n$.

A permutation $p_1p_2\ldots p_n$ {\em contains a patten} $q=q_1\ldots q_k$ 
if $k\leq n$ and some subsequence $p_{i_1}\ldots p_{i_k}$ is order isomorphic to $q$.
A permutation {\em avoids} a pattern if it does not contain it.

Given a set of patterns $q_1,q_2,\ldots$ the set of all permutations which 
avoid all of the patterns is called a {\em pattern-avoiding class}, and the set 
$q_1,q_2,\ldots$ is called its {\em basis}.
A class is {\em finitely based} if there are finitely many $q_i$ defining it.

Given a basis, one would like to know the number of permutations of each length that avoid it. Define
$s_n(B)$ to be the number of permutations of length $n$ which avoid all the patterns in $B$.
Even for bases of single patterns computing this sequence is a difficult problem, 
and results are known only for quite special cases. One can define a formal power 
series or  {\em generating function}
for this sequence.  As mentioned in the introduction, the generating function for $s_n(1234)$ is not algebraic. 
The asymptotic behaviour of the number of 1234-avoiding permutations was                                                                                               
given by Regev \cite{\regev}, and is incompatible with the algebraicity of the associated                                                                                            
generating function \cite{flajolet}. It is however D-finite.
An open question is whether for every pattern $q$, the generating
 function for $s_n(q)$ is D-finite. See \cite{\Stanley} for background on D-finiteness.

\begin{defn}[Linear-bounded automaton]
A {\em \lba} is a nondeterministic Turing machine
such that for every input $w \in \Sigma^*$, there is some accepting
computation in which the tape contains at most $|w|+1$ symbols.
\end{defn}

\begin{defn}[Context-sensitive]
A {\em context-sensitive grammar} consists of three finite sets
$N,T,P$ where $N$ is a set of {\em non-terminals} which are denoted by
upper case letters, and includes a distinguished letter $S$ called the
{\em start symbol}, $T$ is a set of {\em terminals} which are denoted
by lower case letters, and $P$ is a set of {\em productions} of the
form $\alpha\ra \beta$ where $\alpha \in (N\cup T)^*N(N\cup T)^*$ and
$\beta\in (N\cup T)^*$.  The {\em language} of a \cs\ grammar is the set of all
words in $T^*$ that can be obtained by applying some sequence of
productions starting form $S$ until all letters are terminals.
\end{defn}

For readers familiar with \cf\ languages, a context-sensitive
grammar is like a context-free grammar where some productions or rules
can only be applied in the correct ``context''. Just as context-free
languages have a machine-theoretic counterpart (pushdown automata), so
do \csl s.

\begin{lem}\label{lem:lbagrammar}
A \lan\ is generated by a \cs\ grammar if and only if
it is accepted by a \lba.
\end{lem}

 Examples of \csl s are $\{a^nb^nc^n| n\geq 0\}$,
  $\{(ab^n)^n| n\geq 0\}$ and $\{1^p| p$ is prime$\}$ (which are not
  context-free). Context-sensitive languages are closed under union,
  concatenation, star, intersection and complement \cite{\Imm}.

There are many other formal languages of interest, that can be described in terms of grammars or machines of various types. See \cite{\HU}. We breifly mention stack languages since they have some relevance to this paper.

A stack automaton consists of a finite state control, a stack and a pointer to the stack, such that  one can add and delete tokens from the top of the stack, as well as move the pointer down into the stack to read-only. There is a finite set of input letters, stack tokens, and finitely many transition rules. The transisitions are given by a function taking a state of the finite state control and an input letter, and outputing a stack instruction and a new state. If the input is read keeping the stack operational and the final state is an accept state, then we accept, otherwise reject.

\subsection{Non-D-finite stack and \csl s }

As noted in the introduction, the fact that regular and unambiguous \cf\ languages have rational and algebraic generating functions respectively, leads one to ask whether this connection continues for other classes of languages and generating functions. Here we give two examples of \csl s whose generating functions are not D-finite. The question of whether there are ambiguous \cf\ languages with non-D-finite or even non-algebraic generating functions, to the authors knowledge at least, is open.

Consider the language
$\{a^{i_1}b^{i_2}a^{i_3}\ldots a^{i_{2k}}(b^{i_{2k+1}})   \; | \; i_1\leq i_2 \leq i_3 \ldots \leq i_{2k+1}\}$. While this language is not \cf, is can be recognised using a stack automaton as follows. 

Place a $\$$ token on the bottom of the stack. Read $a$ and place a $1$ token on the stack. Repeat until the first $b$ is read. For each $b$, move the pointer down the stack one token. If an $a$ is read before the pointer is at the bottom of the stack, reject. Otherwise if the bottom of the stack is reached, return the pointer to the top of the stack, and for each $b$, place a $1$ token on the stack. Repeat reversing the roles of $a$ and $b$, and if the end of the word is reached then accept.

It follows that this language is indexed, and hence \cs. 
The number of strings in this language of length $n$ is equal to the
number of ways one can partition $n$ as a sum of smaller integers. Since the generating function for this set has infinitely many singularities, it cannot be  D-finite \cite{\flajolet}.
This language is similar to that in Girgorchuk and Machi \cite{\GrigIndex}.

Next, consider the language $\{a^n\;|\; n $ is prime $\}$. This language can be recognized by a \lba\ as follows. Start with the word $aa\ldots a$ written on the tape, with markers as either end. We need to test if the number of $a$s has any divisors. Place a mark $\ast$ on the second $a$. This gives you a way to count ``$2$" on the tape. Suppose the marked $a^{\ast}$ is $i$ squares from the start. Use this to mark every $i$-th square after $a^{\ast}$. If you reach the end of the tape and have a multiple of $i$ $a$'s, then reject, else move the mark $\ast$ up one square to count $i+1$. (This is effectively the sieve method).

The fact that this language is \cs\ is well known and a more careful proof can be found in the literature.
The number of strings in this language of length $n$ is $1$ is $n$ is prime and $0$ otherwise. The generating function for such a set will not satisfy any polynomial recurrence relation, since there are arbitrarily long gaps between successive primes, so this language is not P-recursive so not D-finite.

The fact that there are arbitrarily long gaps between primes follows from the prime number theorem, but a cute argument is that 
 the sequence $k!+2,k!+3,\ldots, k!+k$ is a gap of length $k-1$ for any $k$, since each term has a divisor of $2,3,\ldots,k$ respectively.

\section{The encoding}\label{sec:encoding}

We begin this section with a description of the ``insertion encoding''
 \cite{\ALR}. The idea of this encoding is to build up a \perm\ by
 successively inserting the next highest entry in some open {\em
 slot}, starting from a single open slot, until all slots are filled.
 One can insert in one of four different ways: one can insert on the
 {\em left} of a slot, on the {\em right}, in the {\em middle}
 (creating two slots from one) or {\em filling} the slot.

As an example, the instructions ``middle,right,left,fill,fill'' or
$mrlff$ build the \perm\ $34215$ as follows:
$$\diamond
\ra  \diamond \; 1 \; \diamond
\ra \diamond\; 2\;  1 \; \diamond
\ra  3 \; \diamond\; 2\;  1 \; \diamond
\ra  3 \; 4\; 2\;  1 \; \diamond
\ra  3 \; 4\; 2\;  1 \; 5
$$
In this example, not knowing any better, we performed each insertion
on the left-most open slot. If we want to use any other slot, we
precede the instruction by some number of letter $t$ for {\em
translate}, such that you shift one slot to the right for each $t$
preceding each instruction.

So for the above example, the instructions $mrtltff$ build the
\perm\ $52134$ as follows:
$$\diamond
\ra  \diamond \; 1 \; \diamond
\ra  \diamond\; 2\;  1 \; \diamond
\ra   \diamond\; 2\;  1 \; 3\; \diamond
\ra   \diamond\; 2\;  1 \; 3\; 4
\ra   5\; 2\;  1 \; 3\; 4
$$

The original insertion encoding \cite{\ALR} used integer subscripts on
the $l,r,m,f$ letters to indicate the choice of slot. This means that the 
alphabet is not finite in general. For patterns such as 312 they proved that one
can only ever insert in the first slot to avoid the pattern, so in this case only the four 
(unsubscripted) letters are required. 

We call the
present version the {\em $it$-encoding} which stands for ``insert and
translate''. 

\begin{defn}[$it$-encoding]
Define a codeword be a string of letters $l,r,m,f,t$, subject to the
condition that the number of $t$ letters immediately preceding any of
the other letters $x\in \{l,r,m,f\}$ is at most $k$ where $k$ is the
difference between the number of $m$ and the number of $f$ letters
that appear to the left of $x$, and the codeword ends with an
$f$. Then such a codeword represents a \perm\ $p$ obtained by a
sequence of insertions defined by the following procedure.

Define the {\em next entry} to be the number that is to be inserted
next into the \perm\ and denote it by $i$, define the {\em next
slot} to be number of the slot (numbered from left to right) in $p$
into which the insertion will take place, and the {\em next letter} to
be the next letter of the codeword.

Set  $p=\diamond$, $i=1$ and the next slot to be $1$.

Repeat until the last letter has been read:

Read the next letter of the codeword. 

\bi \item If the next letter is $t$, increment the next slot by $1$.
\item If the next letter is $l$ then replace the next slot by
$i \diamond$. Increment $i$ by $1$ and reset the next slot back
to $1$.
\item If the next letter is $r$ then replace the next slot by
$\diamond i$. Increment $i$ by $1$ and reset the next slot back
to $1$.
\item If the next letter is $m$ then replace the next slot by $\diamond
i\diamond$ (thus increasing the number of slots in $p$ by 1).  Increment
$i$ by $1$ and reset the next slot back to $1$.
\item If the next letter is $f$ then replace the next slot by $i$
(thus decreasing the number of slots in $p$ by 1). Increment $i$ by
$1$ and reset the next slot back to $1$.  \ei
\end{defn}

Note that the length of the \perm\ created is equal to the number of
non-$t$ letters in the code-word. Thus to enumerate permutations in
this encoding of a given length, we would make use of a two-variable
generating function, with one variable for non-$t$ letters and the
other for the $t$s. Then we could set the second variable to $1$.

\section{$E(q)$ is context-sensitive}\label{sec:main}

Let denote $E$  the set of all legal codewords and $E(q)$  the set of codewords representing
permutations that {\em avoid} the pattern $q$. 
In this section we prove  that $E(q)$ is \cs.
 The first step is to show that a
\lba\ can check that a string of $l,r,m,f,t$ is a legal codeword.

\begin{lem}\label{lem:legalcodeword}
Let $w\in\{l,r,m,f,t\}^*$ written on $|w|$ squares of tape and let
$n\in \N$. Then one can determine whether $w$ is a legal code-word
(that is, a word in $E$) for some \perm\ of length $n$ and return the
original string $w$ on the tape, using only the $|w|$ squares in use
and at most $O(|w|^2)$ steps.
\end{lem}

\noindent  \textit{Proof.}
To to this, we need to make sure that the right number of $t$s appear
before each non-$t$ letter. Scan the tape (from left to right) until
you find an $m$. Move right until either you find another $m$, or an
$f$. (if you find neither, reject). When you find an $m,f$ pair (with
only $r,l,t$ letters in between), overwrite the $m$ and $f$ with $m^*$
and $f^*$.

Now scan from $m^*$ to $f^*$ and mark a single $t$ (that is, replace
it by $t^*$) in front of every $r,l$ and the final $f^*$, if there are
$t$s in front of them.

Continue until you cannot find any more unmarked $m,f$ pairs. If there
remains an $m$, reject. If there remains an $f$ which is not at the
end of the tape, reject. (you should still have the $f^{\#}$ at the
end of the tape). If there remain any unmarked $t$ letters, reject.

If not, then the work has the correct number of $t$s in front of each
non-$t$, and the input ends with an $f$, so represents a correct
\perm.
$\Box$

If $x\in \{l,r,m,f\}$ occurs in a code-word $w$ then the entry of the
encoded \perm\ $p$ corresponding to $w$ inserted by $x$ is
denoted by $p_x$.

\begin{lem}\label{lem:pxlesspy}
If $w=w_1xw_2yw_3$ is a codeword with $x,y\in\{l,r,m,f\}$ then
$p_x<p_y$.
\end{lem}

\noindent  \textit{Proof.}
Since entries are inserted in increasing order, and $x$ occurs before
$y$ in the codeword, then the entry inserted by $x$ is smaller than
the entry that is inserted at a later time by $y$.
$\Box$

We will make use of the following lemma.
\begin{lem}\label{lem:xy}
Let $w$ be an $it$-encoded word written on $|w|$ squares of tape and
let $x,y\in \{l,r,m,f\}$ be such that $w=w_1xw_2yw_3$. Then one can
determine whether $p_x$ occurs before or after $p_y$ in the \perm\
corresponding to $w$, (that is, determine whether the sub-permutation
consisting of $p_x$ and $p_y$ is order isomorphic to $12$ or $21$),
and return the original string $w$ on the tape, using only the $|w|$
squares the input is written on and at most $O(|w|^2)$ steps.
\end{lem}

\noindent  \textit{Proof.}
By the previous lemma, since $x$ is to the left of $y$, $p_x<p_y$. So
$x$ inserts the ``1'' and $y$ inserts the ``2''.

The slot in which $x$ is inserted is determined by ``counting'' the
preceding $t$s, and we can do this by marking each preceding $t$ by a
$^*$. Suppose there are $i$ marked $t$s, which means $x$ inserts into
the $(i+1)$-th slot. If $x$ is an $l$ or $f$ then there are $i$ slots
to the left of $x$. If $x$ is an $r$ or $m$ then there are $(i+1)$
slots to the left of $x$, so star $x$ as well.

So the number of starred squares will represent the number of open
slots to the left of $x$. This number can be changed if an $m$ or $f$
inserts in a slot to the left of $x$.

 Move right towards $y$. For each $m$ passed, if we determine that $m$
inserts into a slot before the position of $x$, we need to up our
count by $+1$, (that is, mark $m$ by a $^*$ to signify a shift of
slots). To check this, compare the starred entries with the $t$s
immediately before this $m$. At the machine level, this can be
achieved by scanning to the first $t$ to the left of this $m$ and
marking it by $^{\dag}$, then scan left to the first starred square,
which you remark with $^{**}$, then move right to the $t$ to the left
of $t^{\dag}$, then left and so on until either you run out of starred
squares or run out of $t$s preceding the $m$. If there are more
starred entries, then this $m$ indeed occurs before $x$, so the
position of $x$ needs to be adjusted (so just star $m$). If $m$ has
more, do nothing. Unmark the $t^{\dag}$s and remark the $^{**}$
squares with $^*$.

If you encounter an $f$ before you get to $y$, again if it is
determined that this $f$ fills a slot before the position of $x$ then
we remove one of the stars from the rightmost starred square. We can
determine this by comparing starred squares with $t$s preceding this
$f$ by the same procedure used in the previous paragraph.

If at any stage there are no starred entries left, then $x$ now lies
to the left of any open slots, so $p_x,p_y$ give a $12$
sub-permutation.

So by the time you get to $y$, the position of $x$ is preserved.  Now
compare the number of $t$s preceding $y$ with the number of starred
entries. If there are $j$ $t$s preceding $y$, and $i$ starred
squares, then $y$ inserts into the $(j+1)$-th slot, so if $i\geq j+1$
then $p_x,p_y$ give a $21$ sub-permutation, and if $i<j+1$ then
$p_x,p_y$ give a $12$ sub-permutation. Note that the type of $y$
($l,r,m,f$) doesn't affect this step.

Each time we pass an $m$ or $f$ and when we reach $y$, the machine
scans left and right (marking with $^{**}$ and $^{\dag}$) at most
$|w|$ times, and there are no more than $|w|$ $m,f$s between $x$ and
$y$, so this procedure takes $O(|w|^2)$ steps.
$\Box$

\begin{thm}\label{thm:main} 
$E(q)$ is accepted by a linear-bounded automaton in polynomial time.
\end{thm}

\noindent  \textit{Proof.}
Let $|q|=k$.  We will describe the \lba\ that can decide whether to
accept or reject a string of $l,r,m,f,t$ letters written on the tape,
using the tape alphabet $l,r,m,f,t,l^*,r^*,m^*,f^*,t^*$.

First, verify that the input is a legal code-word by enacting the
procedure given in Lemma \ref{lem:legalcodeword}. Note that the last
letter must be an $f$. This takes at most $O(|w|^2)$ steps.

Next, we systematically run through all ordered $k$-tuples of squares,
such that each square selected has a non-$t$ letter on it. Let
$\overline{x}=(x_1,x_2\ldots, x_k)$ be such a tuple.

We want to test if $\overline{x}$ is order isomorphic to $q$. We do
this by performing pairwise comparisons using the procedure of Lemma
\ref{lem:pxlesspy}. This involves $k\!$ comparisons, but since $q$ is
a fixed permutation, this is a constant.

Each comparison takes $O(|w|^2)$ time and we must repeat this
${|w|\choose k}\in O(|w|^k)$ times, so the algorithm described here
runs in polynomial time.

If a sub-permutation that is order isomorphic to $q$ is found, reject,
else accept.
$\Box$

Since \csl s are closed under intersection \cite{\HU} then any finitely based class 
is \cs. We prove this fact explicitly.
\begin{cor}
The set of permutations which avoid a finite list of forbidden subpatterns 
can be encoded as a language that is 
accepted by a linear-bounded automaton in polynomial time.
\end{cor}

\noindent  \textit{Proof.}
Suppose $q_1,\ldots, q_k$ are a finite list of forbidden subpatterns. By Theorem \ref{thm:main}
we can construct linear-bounded automata to accept each of $E(q_i)$, such that 
each automaton leaves the original 
codeword on the tape, and runs in polynomial time. Start with a codeword written on the tape, and run 
the automaton to decide $E(q_1)$. The automaton accepts or rejects, and leaves the codeword on the tape. If
it accepts, run the automaton to decide $E(q_2)$ and so on. If each of the $k$ automata accept, then accept, else reject.
The time taken is a constant ($k$) times a polynomial.
$\Box$

\subsection{Full encoding  is accepted by a stack automaton}
In this section we will show that $E$ is accepted by a stack
automaton.

\begin{prop}
$E$ is accepted by a stack automaton.
\end{prop}
\noindent  \textit{Proof.}
 We will describe a deterministic stack automaton that accepts the
set of all strings of $l,r,m,f,t$ that are codewords. That is, we must
check that the number of $t$s immediately preceding any of the other
letters is no more that the difference between the number of $m$s and
$f$s preceding it.

Recall that a stack automaton constists of a finite state control, and single stack, and a pointer which 
points to some entry in the stack.
 There are three main states, $q_0,q_A,q_F$, start,
accept and fail respectively. 
Place the entry {\em root} on the bottom of the stack.

First we describe the transition function for the start state
$q_0$. Let $x$ be any entry of the stack, and let $x_r$ be any entry
except the root.  

\bi \item
$\delta(q_0,l,x)$: move the pointer to the top of
the stack (that is, go into a state $q_{up}$ and perform single
$\epsilon$ transitions moving the cursor up one step) and continue.

\item
$\delta(q_0,r,x)$: move the pointer to the top of the stack and
continue.

\item $\delta(q_0,m,x)$: move the pointer to the top of the stack, add
a new entry, and continue (so the pointer points to the
new entry).

\item $\delta(q_0,f,x_r)$: move the pointer to the top of the stack,
delete the top entry, (so the pointer reverts back one step) and
continue.

\item $\delta(q_0,f,\mathit{root})$: if all the input has been read,
move to $q_A$, or if not, move to $q_F$, and stop.

\item $\delta(q_0,t,x_r)$: move the pointer down the stack one step.

\item $\delta(q_0,t,\mathit{root})$: move to $q_F$, and stop.  \ei

The machine reads input letters from left to right, and the height of
the stack indicates the current sum of number of $m$s minus
$f$s. When a $t$ is read, the machine checks that this sum is not
zero (the pointer points above the root). For each $t$ read in
succession, the pointer moves down the stack, so if the number of
$t$s preceding another letter exceeds $(m-f)$ the input is rejected,
otherwise, if the end of the input is reaches, and the last letter is
an $f$ with an empty stack, accept. If an $f$ is read on an empty
stack before all the input has been read, then no slots remain so the
string is not a codeword.
$\Box$

\section{Acknowledgements}

Thanks to Mireille Bousquet-Melou,  Marni Mishna,
Andrew Rechnitzer, Mike Zabrocki and also to an anonymous reviewer 
for their careful reading of an earlier draft, and their invaluable 
feedback and ideas for this paper.

\bibliography{refs} \bibliographystyle{elsart-num}

\end{document}